\documentclass[a4,twocolumn]{article}

\usepackage{graphicx}

\def\be{\begin{equation}}
\def\ee{\end{equation}}
\newcommand{\kk}[2]{\frac{#1}{#2}}

\newcommand{\ff}[1]{{\bf #1}}
\def\a{\alpha}

\def\e{\epsilon}
\def\lam{\lambda}
\def\s{\qquad}

\def\={\approx}
\def\8{{\infty}}
\def\x{\ff{x}}
\def\vcode#1#2#3#4{\begin{figure}
\begin{center}
\begin{minipage}[c]{#1\textwidth}
{{\small #2 \hrule \vspace{5pt} {\bf begin}  \\   %
{\it #3} \\ {\bf end} \vspace{5pt} \hrule }}
\end{minipage}
\caption{#4}
\end{center}   \end{figure}} 

\begin{document}

\title{Cuckoo Search via L\'evy Flights}

\author{Xin-She Yang \\
Department of Engineering, \\
University of Cambridge \\
Trumpinton Street \\
Cambridge CB2 1PZ, UK \\
\and
Suash Deb \\
Department of Computer Science \& Engineering \\
C. V. Raman College of Engineering \\
Bidyanagar, Mahura, Janla \\
Bhubaneswar 752054, INDIA \\
}

\date{{\bf Citation detail}: X.-S. Yang, S. Deb,
``Cuckoo search via L\'evy flights'', in: {\it Proc. of \\
World Congress on Nature \& Biologically Inspired Computing} (NaBIC 2009), \\
December 2009, India. IEEE Publications, USA,  pp. 210-214 (2009).  }

\maketitle

\abstract{

In this paper, we intend to formulate a new metaheuristic algorithm, called Cuckoo Search (CS),
for solving optimization problems. This algorithm is based on
the obligate brood parasitic behaviour of some cuckoo species in combination with the L\'evy flight
behaviour of some birds and fruit flies. We validate
the proposed algorithm against test functions and then compare its performance with
those of genetic algorithms and particle swarm optimization. Finally, we discuss the
implication of the results and suggestion for further research. \\

}



\section{Introduction}

More and more modern metaheuristic algorithms inspired by nature are emerging
and they become increasingly popular. For example, particles swarm optimization (PSO)
was inspired by fish and bird swarm intelligence, while the Firefly Algorithm
was inspired by the flashing pattern of tropical fireflies \cite{Bod,Blum,Deb,Yang,Yang2}.
These nature-inspired metaheuristic algorithms have been used in a wide range
of optimization problems, including NP-hard problems such
as the travelling salesman problem \cite{Bod,Blum,Deb,Gold,Ken2,Yang}.

The power of almost all modern metaheuristics comes from the fact that they
imitate the best feature in nature, especially biological systems
evolved from natural selection over millions of years. Two important characteristics
are selection of the fittest and adaptation to the environment. Numerically speaking,
these can be translated into two crucial characteristics of the modern
metaheuristics: intensification and diversification \cite{Blum}. Intensification
intends to search around the current best solutions and select the best
candidates or solutions, while diversification makes
sure the algorithm can explore the search space efficiently.

This paper aims to formulate  a new algorithm, called Cuckoo Search (CS),
based on the interesting breeding bebaviour such as brood parasitism
of certain species of cuckoos.
We will first introduce the breeding bebaviour of cuckoos and the characteristics of
L\'evy flights of some birds and fruit flies, and then formulate the
new CS, followed by its implementation. Finally, we will compare the
proposed search strategy with other popular optimization algorithms and
discuss our findings and their implications for various optimization problems.

\section{Cuckoo Behaviour and L\'evy Flights}

\subsection{Cuckoo Breeding Behaviour}

Cuckoo are fascinating birds, not only because of the beautiful sounds they can make,
but also because of their aggressive reproduction strategy. Some species such as the
{\it ani} and {\it Guira} cuckoos lay their eggs in communal nests, though they may remove others'
eggs to increase the hatching probability of their own eggs \cite{Payne}. Quite a
number of species engage the obligate brood parasitism by laying their eggs in the
nests of other host birds (often other species). There are three basic types of
brood parasitism: intraspecific brood parasitism, cooperative breeding, and
nest takeover. Some host birds can engage direct conflict with the intruding
cuckoos. If a host bird discovers
the eggs are not their owns, they will either throw these alien eggs away
or simply abandon its nest and build a new nest elsewhere. Some cuckoo species
such as the New World brood-parasitic {\it Tapera}
have evolved in such a way that female parasitic cuckoos are often
very specialized in the mimicry in colour and pattern of the eggs of a few chosen
host species \cite{Payne}. This reduces the probability of their eggs
being abandoned and thus increases their reproductivity.

In addition, the timing of egg-laying of some species
is also amazing. Parasitic cuckoos often choose a nest where the
host bird just laid its own eggs. In general,
the cuckoo eggs hatch slightly earlier than their host eggs. Once
the first cuckoo chick is hatched, the first instinct action it
will take is to evict the host eggs by blindly propelling
the eggs out of the nest, which increases the cuckoo chick's
share of food provided by its host bird. Studies also show that a cuckoo chick can
also mimic the call of host chicks to gain access to more feeding opportunity.

\subsection{L\'evy Flights}

On the other hand, various studies have shown that
flight behaviour of many animals and insects has demonstrated
the typical characteristics of L\'evy flights
\cite{Brown,Reynolds,Pav,Pav2}. A recent study by Reynolds and Frye shows that
fruit flies or {\it Drosophila melanogaster}, explore their landscape using a series
of straight flight paths punctuated by a sudden $90^{o}$ turn, leading to
a L\'evy-flight-style intermittent scale free search pattern.
Studies on human behaviour such as the Ju/'hoansi hunter-gatherer foraging
patterns also show the typical feature of L\'evy flights.
Even light can be related to L\'evy flights \cite{Barth}.
Subsequently, such behaviour has been applied to
optimization and optimal search, and preliminary results show its
promising capability \cite{Pav,Reynolds,Shles,Shles2}.

\section{Cuckoo Search}

For simplicity in
describing our new Cuckoo Search, we now use the following
three idealized rules: 1) Each cuckoo lays one egg at a time, and dump
its egg in randomly chosen nest; 2) The best nests with high quality of
eggs will carry over to the next generations;
3) The number of available host nests
is fixed, and the egg laid by a cuckoo is discovered by the host bird
with a probability $p_a \in [0,1]$. In this case, the host bird can either
throw the egg  away or abandon the nest, and build a completely new nest. For simplicity,
this last assumption can be approximated by the fraction $p_a$ of the $n$ nests
are replaced by new nests (with new random solutions).

For a  maximization problem, the quality or fitness of a solution can simply be proportional
to the value of the objective function. Other forms of fitness can be defined in a similar
way to the fitness function in genetic algorithms. For simplicity, we can use the following
simple representations that each egg in a nest represents a solution, and a cuckoo egg
represent a new solution, the aim is to use the new and potentially better solutions (cuckoos)
to replace  a not-so-good solution in the nests. Of course, this algorithm can be extended to
the more complicated case where each nest has multiple eggs representing
a set of solutions. For this present work, we will use the simplest approach where each nest has only
a single egg.

Based on these three rules, the basic steps of the Cuckoo Search (CS)
can be summarized as the pseudo code shown in Fig. \ref{kuk-fig-100}.

\vcode{0.45}{{\sf Cuckoo Search via L\'evy Flights}} {
\indent \quad Objective function $f(\x), \;\; \x=(x_1, ..., x_d)^T$ \\
\indent \quad Generate initial population of \\
\indent \qquad \quad $n$ host nests $\x_i \; (i=1,2,...,n)$ \\
\indent \quad {\bf while} ($t<$MaxGeneration) or (stop criterion) \\
\indent \qquad Get a cuckoo randomly by L\'evy flights \\
\indent \qquad \quad evaluate its quality/fitness $F_i$  \\
\indent \qquad Choose a nest among $n$ (say, $j$) randomly \\
\indent \qquad  {\bf if } ($F_i>F_j$), \\
\indent \qquad \qquad replace $j$ by the new solution;\\
\indent \qquad  {\bf end} \\
\indent \qquad  A fraction ($p_a$) of worse nests \\
\indent \qquad \qquad are abandoned and new ones are built; \\
\indent \qquad Keep the best solutions \\
\indent \qquad \qquad (or nests with quality solutions); \\
\indent \qquad Rank the solutions and find the current best  \\
\indent \quad {\bf end while} \\
\indent \quad Postprocess results and visualization }{Pseudo code of
the Cuckoo Search (CS). \label{kuk-fig-100} }

When generating new solutions $\x^{(t+1)}$ for, say, a cuckoo $i$,
a L\'evy flight is performed
\be \x^{(t+1)}_i=\x_i^{(t)} + \a \oplus \textrm{L\'evy}(\lam), \ee
where $\a>0$ is the step size which should be related to the scales of the problem of
interests. In most cases,  we can use $\a=1$. The above equation is essentially the
stochastic equation for random walk. In general, a random walk is a Markov chain
whose next status/location only depends on the current location (the first term in the above equation)
and the transition probability (the second term). The product $\oplus$ means entrywise multiplications.
This entrywise product is similar to those used in PSO, but here the random walk
via L\'evy flight is more efficient in exploring the search space
as its step length is much longer in the long run.

The L\'evy flight  essentially provides a random walk while the random step length
is drawn from a L\'evy distribution
\be \textrm{L\'evy} \sim u = t^{-\lam}, \s (1 < \lam \le 3), \ee
which has an infinite variance with an infinite mean.
Here the steps essentially form
a random walk process with a power-law step-length distribution with a heavy tail.
Some of the new solutions should be generated by L\'evy walk around the
best solution obtained so far, this will speed up the local search. However,
a substantial fraction of the new solutions should be generated by far field randomization
and whose locations should be far enough from the current best solution, this
will make sure the system will not be trapped in a local optimum.

From a quick look, it seems that there is some similarity between CS and hill-climbing
in combination with some large scale randomization. But there are some significant differences.
Firstly, CS is a population-based algorithm, in a way similar to GA and PSO, but it uses
some sort of elitism and/or selection similar to that used in harmony search.
Secondly, the randomization
is more efficient as the step length is heavy-tailed, and any large step is possible.
Thirdly, the number of parameters to be tuned is less than GA and PSo, and thus it is
potentially more generic to adapt to a wider class of optimization problems.
In addition, each nest can represent a set of solutions, CS can thus be extended
to the type of meta-population algorithm.

\section{Implementation and Numerical Experiments}

\subsection{Validation and Parameter Studies}

\begin{figure}
 \centerline{\includegraphics[height=2.5in,width=3in]{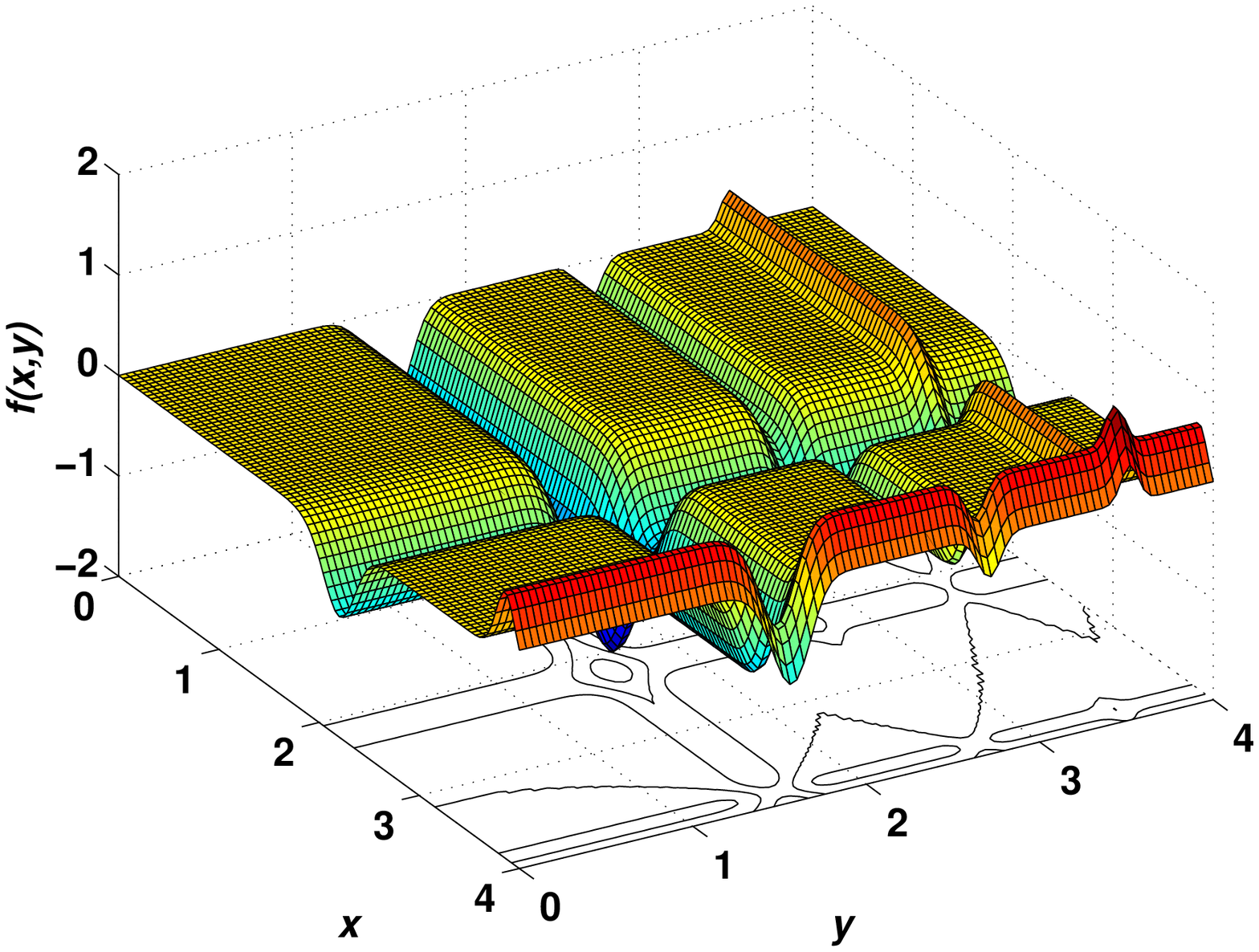} }
\vspace{-5pt}
\caption{The landscaped of Michaelwicz's function. \label{kuk-fig-200} }
\end{figure}

\begin{figure}
 \centerline{\includegraphics[height=2.5in,width=3in]{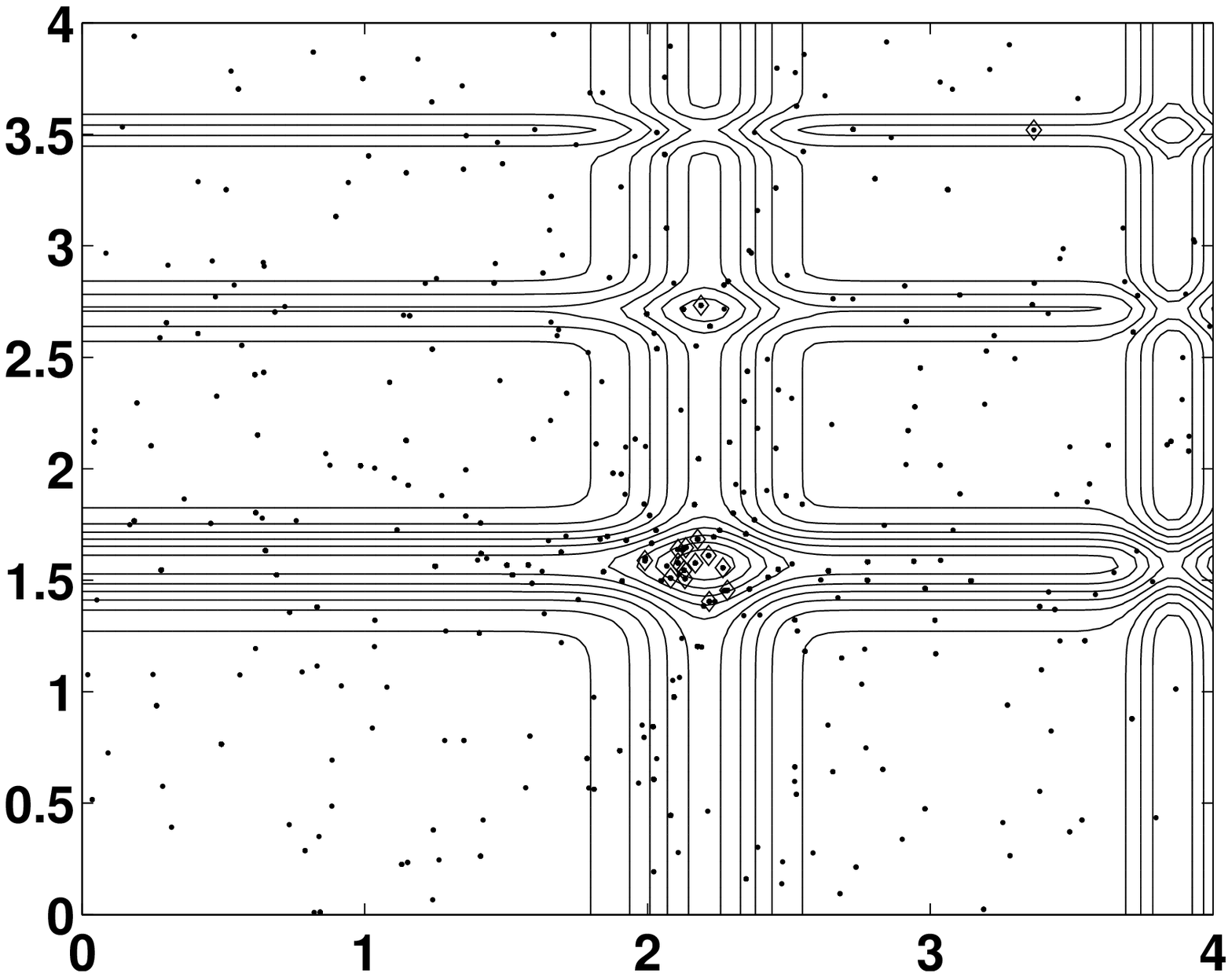}}
\vspace{-5pt}
\caption{Search paths of nests using Cuckoo Search.
The final locations of the nests are marked with {\large $\diamond\!\!\!\;\!\cdot$} in the figure. \label{kuk-fig-300} }
\end{figure}

After implementation, we have to validate the algorithm using test functions with analytical
or known solutions. For example, one of the many test functions we have used is
the bivariate Michaelwicz function
\be f(x,y)=-\sin(x) \sin^{2m}(\kk{x^2}{\pi})-\sin (y) \sin^{2m} (\kk{2 y^2}{\pi}), \ee
where $m=10$ and $(x,y) \in [0,5] \times [0,5]$. This function has a global
minimum $f_* \=-1.8013$ at $(2.20319,1.57049)$. The landscape of this funciton
is shown in Fig. \ref{kuk-fig-200}. This global optimum can easily be
found using Cuckoo Search, and the results are shown in Fig. \ref{kuk-fig-300}
where the final locations of the nests are also marked with {\large $\diamond\!\!\cdot$ }
in the figure. Here we have
used $n=15$ nests, $\a=1$ and $p_a=0.25$.  In most of our simulations, we have used
$n=15$ to $50$.

From the figure, we can see
that, as the optimum is approaching, most nests aggregate towards the global optimum.
We also notice that the nests are also distributed at different (local) optima in the
case of multimodal functions. This means that CS can find all the optima simultaneously
if the number of nests are much higher than the number of local optima. This advantage
may become more significant when dealing with multimodal and multiobjective optimization
problems.

We have also tried to vary the number of host nests (or the population size $n$)
and the probability $p_a$. We have used $n=5,10,15,20,50, 100,150, 250, 500$
and $p_a=0, 0.01, 0.05, 0.1, 0.15, 0.2, 0.25$, $0.4, 0.5$.
From our simulations, we found that $n=15$ and $p_a=0.25$ are sufficient for
most optimization problems. Results and analysis also imply that the convergence rate, to
some extent, is not sensitive to the parameters used. This means that the fine adjustment
is not needed for any given problems.
Therefore, we will use fixed $n=15$ and $p_a=0.25$ in the
rest of the simulations, especially for the comparison studies presented in the
next section.

\subsection{Test Functions}

There are many benchmark test functions in literature \cite{Chatt,Shang,Schoen}, and they are designed to
test the performance of optimization algorithms. Any new optimization algorithm should also be
validated and tested against these benchmark functions.  In our simulations, we have used the
following test functions.

De Jong's first function is essentially a sphere function
\be f(\x)=\sum_{i=1}^d x_i^2, \s x_i \in [-5.12,5.12], \ee
whose global minimum $f_*=0$ occurs at $\x_*=(0,0,...,0)$.
Here $d$ is the dimension.

Easom's test function is unimodal
\be f(x,y)=-\cos(x) \cos(y) \exp[-(x-\pi)^2 - (y-\pi)^2], \ee
where \be (x,y) \in [-100,100] \times [-100,100]. \ee
It has a global minimum of $f_*=-1$ at $(\pi,\pi)$ in a very small region.

Shubert's bivariate function
\be f(x,y)=-\sum_{i=1}^5 i \cos [(i+1) x +1] \sum_{i=1}^5 \cos [(i+1) y +1)], \ee
has 18 global minima in the region $(x,y) \in [-10,10] \times [-10,10]$. The value of its
global minima is $f_*=-186.7309$.

Griewangk's test function has many local minima
\be f(\x)=\kk{1}{4000} \sum_{i=1}^d x_i^2 - \prod_{i=1}^d \cos(\kk{x_i}{\sqrt{i}})+1, \ee
but a single global mimimum $f_*=0$ at $(0,0,...,0)$ for all $-600 \le x_i \le 600$
where $i=1,2,...,d$.

Ackley's function is multimodal
\[ f(\x)=-20 \exp\Bigg[-0.2 \sqrt{\kk{1}{d} \sum_{i=1}^d x_i^2} \; \Bigg] \]
\be -\exp[\kk{1}{d} \sum_{i=1}^d \cos (2 \pi x_i)]
+ (20+e), \ee
with a global minimum $f_*=0$ at $\x_*=(0,0,...,0)$ in the
range of $-32.768 \le x_i \le 32.768$ where $i=1,2,...,d$.

The generalized Rosenbrock's function is given by
\be f(\x)=\sum_{i=1}^{d-1} \Big[ (1-x_i)^2 + 100 (x_{i+1}-x_i^2)^2 \Big], \ee
which has a minimum $f(\x_*)=0$ at $\x_*=(1,1,...,1)$.

Schwefel's test function is also multimodal
\be f(\x)=\sum_{i=1}^d \Big[ -x_i \sin (\sqrt{|x_i|} \;) \Big],
\;\; -500 \le x_i \le 500, \ee
with a global minimum of $f_*=-418.9829d$ at $x_i^*=420.9687 (i=1,2,...,d)$.

Rastrigin's test function
\be f(\x)=10 d + \sum_{i=1}^d [x_i^2 - 10 \cos (2 \pi x_i) ], \ee
has a global minimum $f_*=0$ at $(0,0,...,0)$ in a hypercube
$-5.12 \le x_i \le 5.12$ where $i=1,2,...,d$.

Michalewicz's test function has $d!$ local optima
\be f(\x)=-\sum_{i=1}^d \sin (x_i) \Big[\sin (\kk{i x_i^2}{\pi}) \Big]^{2m}, \;\; (m=10), \ee
where $0 \le x_i \le \pi$ and $i=1,2,...,d$. The global mimimum is
$f_*\=-1.801$ for $d=2$, while $f_* \=-4.6877$ for $d=5$.

\subsection{Comparison of CS with PSO and GA}

Recent studies indicate that PSO algorithms can
outperform genetic algorithms (GA) \cite{Gold}
and other conventional algorithms for many optimization problems.
This can partly be attributed to the broadcasting ability of the current
best estimates which potentially gives better and quicker convergence towards the
optimality.  A general framework for evaluating statistical
performance of evolutionary algorithms
has been discussed in detail by Shilane et al. \cite{Shilane}.

Now we will compare the Cuckoo Search with PSO and genetic
algorithms for various standard test functions.
After implementing these algorithms using
Matlab, we have carried out extensive simulations and each algorithm has been
run at least 100 times so as to carry out meaningful statistical analysis.
The algorithms stop when the variations of function values are less than
a given tolerance $\e \le 10^{-5}$. The results are
summarized in the following tables (see Tables 1 and 2) where the global optima
are reached. The numbers are in the
format: average number of evaluations (success rate), so
$927 \pm 105 (100\%)$ means that the average number (mean) of function
evaluations is 927 with a standard deviation of 105. The success rate
of finding the global optima for this algorithm is $100\%$.

\begin{table}[ht]
\caption{Comparison of CS with genetic algorithms}
\centering
\begin{tabular}{ccccc}
\hline \hline
Functions/Algorithms & GA  & CS \\
\hline
 Multiple peaks & $52124 \pm 3277 (98\%)$   & $927 \pm 105 (100 \%)$ \\

 Michalewicz's ($d\!\!=\!\!16$)  & $89325 \pm 7914 (95 \%)$  & $3221 \pm 519 (100\%)$ \\

 Rosenbrock's ($d\!\!=\!\!16$) & $55723 \pm 8901 (90\%)$ & $5923 \pm 1937 (100\%) $ \\

 De Jong's ($d\!\!=\!\!256$) & $25412 \pm 1237 (100\%)$  & $4971 \pm 754 (100\%)$\\
 Schwefel's ($d\!\!=\!\!128$) & $227329 \pm 7572 (95\%)$  & $8829 \pm 625 (100\%)$ \\

 Ackley's ($d\!\!=\!\!128$) & $32720 \pm 3327 (90\%)$ & $4936 \pm 903 (100\%)$ \\

 Rastrigin's & $110523 \pm 5199 (77 \%)$ & $10354 \pm 3755 (100\%)$ \\

 Easom's & $19239 \pm 3307 (92\%)$  & $6751 \pm 1902 (100\%)$ \\

 Griewank's & $70925 \pm 7652 (90\%)$  & $10912 \pm 4050 (100\%)$ \\

 Shubert's (18 minima) & $54077 \pm 4997 (89\%)$  & $9770 \pm 3592 (100\%)$ \\

\hline
\end{tabular}
\end{table}

We can see that the CS is much more efficient in finding the global optima
with higher success rates. Each function evaluation is virtually instantaneous
on modern personal computer. For example, the computing time for 10,000 evaluations
on a 3GHz desktop is about 5 seconds.

\begin{table}[ht]
\caption{Comparison of CS with Particle Swarm Optimisation}
\centering
\begin{tabular}{ccccc}
\hline \hline
Functions/Algorithms & PSO & CS \\
\hline
 Multiple peaks  & $3719 \pm  205 (97\%)$ & $927 \pm 105 (100 \%)$ \\

 Michalewicz's ($d\!\!=\!\!16$)  & $6922 \pm 537 (98\%)$  & $3221 \pm 519 (100\%)$ \\

 Rosenbrock's ($d\!\!=\!\!16$)  & $32756 \pm 5325 (98\%)$ & $5923 \pm 1937 (100\%) $ \\

 De Jong's ($d\!\!=\!\!256$)  & $17040 \pm 1123 (100\%)$ & $4971 \pm 754 (100\%)$\\
 Schwefel's ($d\!\!=\!\!128$) & $14522 \pm 1275 (97\%)$ & $8829 \pm 625 (100\%)$ \\

 Ackley's ($d\!\!=\!\!128$)  & $23407 \pm 4325 (92\%)$ & $4936 \pm 903 (100\%)$ \\

 Rastrigin's & $79491 \pm 3715 (90\%)$ & $10354 \pm 3755 (100\%)$ \\

 Easom's  & $17273 \pm 2929 (90\%)$ & $6751 \pm 1902 (100\%)$ \\

 Griewank's  & $55970 \pm 4223 (92\%)$ & $10912 \pm 4050 (100\%)$ \\

 Shubert's (18 minima) & $23992 \pm 3755 (92\%)$ & $9770 \pm 3592 (100\%)$ \\

\hline
\end{tabular}
\end{table}

For all the test functions, CS has outperformed both GA and PSO. The primary reasons are two folds:
A fine balance of randomization and intensification, and less number of control parameters.
As for any metaheuristic algorithm, a good balance of intensive local search strategy
and an efficient exploration of the whole search space will usually lead to a more efficient algorithm.
On the other hand, there are only two parameters in this algorithm, the population size $n$,
and $p_a$. Once $n$ is fixed, $p_a$ essentially controls the elitism and the balance of the
randomization and local search. Few parameters make an algorithm less complex and thus
potentially more generic. Such observations deserve more systematic research and further elaboration
in the future work.

\newpage 
$\;$
\vspace{4.25in}

\section{Conclusions}

In this paper, we have formulated a new metaheuristic Cuckoo Search in combination
with L\'evy flights,  based on the breeding strategy of some cuckoo species. The proposed algorithm has
been validated and compared with other algorithms such as genetic algorithms
and particle swarm optimization. Simulations and comparison show that CS is
superior to these existing algorithms for multimodal objective functions. This
is partly due to the fact that there are  fewer parameters to be fine-tuned in
CS than in PSO and genetic algorithms. In fact, apart from the population size
$n$, there is essentially one parameter $p_a$. Furthermore, our simulations also
indicate that the convergence rate is insensitive to the parameter $p_a$.
This also means that we do not have to fine tune these parameters for
a specific problem. Subsequently, CS is more generic and robust for
many optimization problems, comparing with other metaheuristic algorithms.

This potentially powerful optimization strategy can easily be extended
to study multiobjecitve optimization applications with various constraints, even to
NP-hard problems. Further studies can focus  on the
sensitivity and parameter studies and their possible relationships
with the convergence rate of the algorithm. Hybridization with other popular
algorithms such as PSO will also be potentially fruitful.


\end{document}